\documentclass[11pt,leqno]{article} 
\usepackage{graphics}

\usepackage{lipsum}

\let\OLDthebibliography\thebibliography
\renewcommand\thebibliography[1]{
  \OLDthebibliography{#1}
  \setlength{\parskip}{2pt}
  \setlength{\itemsep}{2pt plus 0.3ex}
}

\newtheorem{thm}{Theorem}[section]
\newtheorem{lma}{Lemma}[section]

\newtheorem{prop}{Proposition}[section]

\newcommand{\beqa}{\begin{eqnarray}}
\newcommand{\eeqa}{\end{eqnarray}}

\newcommand{\pf}{\noindent {\bf Proof:} $\s$ }
\newcommand{\epf}{ \hfill$\diamondsuit$ \medskip}

\newcommand{\beq}{\begin{equation}}
\newcommand{\eeq}{\end{equation}}
\newcommand{\lbl}{\label}
\newcommand{\s}{\; \;}

\newcommand{\la}{\lambda}

\newcommand{\ra}{\rightarrow}
\newcommand{\al}{\alpha}

\title{Unbounded solutions of periodic systems}

\author{
Philip Korman   \\ 
Department of Mathematical Sciences \\ 
University of Cincinnati \\ 
Cincinnati Ohio 45221-0025 \\
}

\date{}

\begin{document}

\maketitle
\begin{abstract}
\noindent
This paper deals with various cases of resonance, which is a fundamental concept of science and engineering.
Specifically, we study the connections between periodic and unbounded solutions for several classes of equations and systems. In particular, 
we extend the classical Massera's theorem, dealing with periodic systems of the type
\[
x'=A(t)x+f(t) \,,
\]
and clarify that this theorem deals with a case of resonance. Then we provide instability results for the corresponding semilinear systems, with the linear part at resonance. We also use the solution curves developed in \cite{K2},\cite{K10} to establish the instability results   for pendulum-like equations, and for first-order periodic equations.
 \end{abstract}

\begin{flushleft}
\end{flushleft}

\begin{flushleft}
\end{flushleft}

\section{Introduction}
\setcounter{equation}{0}
\setcounter{thm}{0}
\setcounter{lma}{0}

We consider equations and systems {\em at resonance}. A textbook example of resonance involves an equation like
\[
x''(t)+x(t)=\sin t \,.
\]
All solutions of the corresponding homogeneous equation are bounded, while the periodic forcing term produces an unbounded response. A similar situation occurs for ($t>0$)
\[
x''(t)+x'(t)=1+\sin t \,.
\]
This is also a case of resonance, for which we shall study nonlinear perturbations, involving pendulum-like equations. We shall also deal with the resonance for first order periodic equations and systems, particularly in the context of Massera's theorem. 
\medskip

Consider  a  system with an $n \times n$ $p$-periodic matrix $A(t)$ and a $p$-periodic vector $f(t) \in R^n$
\beq
\lbl{3in}
x'=A(t)x+f(t) \,,
\eeq
and the corresponding homogeneous system
\beq
\lbl{4in}
x'=A(t)x \,.
\eeq
The famous theorem of J.L. Massera  \cite{ma} says: if (\ref{3in}) has a bounded solution ($||x(t)|| \leq c$ uniformly in $t>0$), then it has a $p$-periodic solution.
The original statement of Massera's theorem is very intriguing, but it appears to be not easy to use. Indeed, if one manages to construct an explicit  bounded solution of a periodic system, chances are that solution is already periodic. We shall deal with the following contraposition form of Massera's theorem.
\begin{thm} (Massera \cite{ma})
\lbl{thm:0}
If (\ref{3in}) has no periodic solution, then all of its solutions are unbounded as $t \ra \infty$. Moreover, (\ref{4in}) has  a $p$-periodic solution.
\end{thm}  
This form appears to be more natural. In particular, it becomes clear that Massera's theorem deals with {\em a case of resonance}, and that this theorem admits a natural extension to a rather complete result, with detailed description of the dynamics of (\ref{3in}) that we present.
\medskip

When studying equations with periodic coefficients a natural first step is to investigate the existence of periodic solutions. What could be the second step? Traditionally, one studies the stability of  periodic solutions, see e.g., B.P. Demidovi\v{c} \cite{d}. Motivated by the second Massera's theorem, see e.g., p. 203 in \cite{K} (or R. Ortega \cite{o} for a detailed presentation), G. Seifert \cite{s} and J.M. Alonso and R. Ortega \cite{O} showed that in case periodic solutions are absent, one can prove that all solutions are unbounded for equations at or near resonance, and that solutions tend to infinity in a stronger norm compared with Massera's theorem. We present similar instability results for semilinear perturbations of linear systems:
\[
x'+A(t)x+f(x)=g(t) \,,
\]
based on the Landesman-Lazer condition (rather than on the Lazer-Leach condition used in \cite{s} and \cite{O}).
\medskip

For a class of pendulum-like equations
of the type
\[
x''(t)+\la x'(t)+g(x)=f(t) \,,
\]
with periodic $f(t)$, and for similar first order equations,
we relied on a detailed description of the  curves of periodic solutions developed in P. Korman \cite{K2}, \cite{K10}, to get conditions that are both necessary and sufficient  for the existence of periodic solutions at resonance, and obtained related instability results.

\section{An  extension of Massera's theorem}
\setcounter{equation}{0}
\setcounter{thm}{0}
\setcounter{lma}{0}

We begin with a single equation
\beq
\lbl{1}
x'(t)+a(t)x(t)=f(t) \,,
\eeq
with continuous $p$-periodic functions $a(t)$ and $f(t)$, so that $a(t+p)=a(t)$ and $f(t+p)=f(t)$ for some $p>0$, and all $t$. Write its general  solution as 
\beq
\lbl{1.1}
x(t)=\frac{1}{\mu (t)} c+\frac{1}{\mu (t)} \int_0^t \mu (s) f(s) \, ds,
\eeq
where $\mu (t)=e^{\int _0^t a(s) \, ds}$, and $c$ is an arbitrary constant. This formula shows that the dynamics is simple in case $\int _0^p a(s) \, ds \ne 0$. Then there exists a unique $p$-periodic solution that attracts all other solutions as $t \ra \infty $ if $\int _0^p a(s) \, ds>0$, and as $t \ra -\infty $, in case $\int _0^p a(s) \, ds<0$, see e.g., \cite{K} for the details. More interesting is the case
\beq
\lbl{2}
\int _0^p a(s) \, ds=0 \,,
\eeq
when the corresponding homogeneous equation
\[
x'+a(t)x=0
\]
has  $p$-periodic solutions $x(t)=\frac{c}{\mu (t)}$, where $\mu (t)$ is $p$-periodic. There are two cases. If $\int_0^p \mu (s) f(s) \, ds=0$ then clearly all solutions of (\ref{1}) are $p$-periodic. In case $\int_0^p \mu (s) f(s) \, ds \ne 0$, all solutions are unbounded as $t \ra \pm \infty$ (just consider $x(mp)$ with $m \ra \pm \infty$, and observe that $x(mp)-x(0)=m \int_0^p \mu (s) f(s) \, ds $, by the periodicity of $\mu (t)$ and $f(t)$). So that the condition (\ref{2}) presents  a full fledged {\em case of  resonance}, even though the equation (\ref{1}) is of first order.
\medskip

We consider now a $p$-periodic system
\beq
\lbl{3}
x'=A(t)x+f(t) \,,
\eeq
and the corresponding homogeneous system
\beq
\lbl{4}
x'=A(t)x \,.
\eeq
We assume that the $n \times n$ matrix $A(t)$ and the vector $f(t) \in R^n$ have continuous entries, and $A(t+p)=A(t)$, $f(t+p)=f(t)$ for some $p>0$ and all $t$.
If $X(t)$ is the fundamental solution matrix of (\ref{4}), then the solution of (\ref{4}) satisfying the initial condition $x(0)=x_0$ is
\[
x(t)=X(t)x_0 \,.
\]
For the non-homogeneous system (\ref{3}), the solution   satisfying the initial condition $x(0)=x_0$, and denoted by $x(t,x_0)$, is given by
\beq
\lbl{4.1}
x(t)=X(t)x_0+X(t) \int_0^t X^{-1}(s) f(s) \, ds \,.
\eeq
The homogeneous system (\ref{4}) has a $p$-periodic solution, with  $x(p)=x(0)$, provided that the $n \times n$ system of linear equations
\beq
\lbl{5}
\left(I-X(p)   \right)x_0=0
\eeq
has a non-trivial solution $x_0$. Define the vector 
\beq
\lbl{5.1}
b=X(p) \int_0^p X^{-1}(s) f(s) \, ds \,. 
\eeq
The non-homogeneous system (\ref{3}) has a $p$-periodic solution, with $x(p)=x(0)$, provided that the system
\beq
\lbl{6}
\left(I-X(p)   \right)x_0=b
\eeq
has a solution $x_0$. If (\ref{3}) has no  $p$-periodic solutions, then the system (\ref{6}) has no solution, so that the matrix $I-X(p)$ is singular. Then (\ref{5}) has non-trivial solutions, and  (\ref{4}) has a $p$-periodic solution.
This justifies the extra claim of  Massera's Theorem \ref{thm:0}. 
\medskip

In the theorem below we shall assume that the homogeneous system (\ref{4}) has a $p$-periodic solution. Then the matrix $X(p)$ has an eigenvalue $1$, and the spectral radius of $X(p)$ is $\geq 1$ (recall that the spectral radius $\rho (X(p))=\max |\la _i| $, maximum taken over all eigenvalues of $X(p)$).

\begin{thm}\lbl{thm:1}
Assume that the homogeneous system (\ref{4}) has a $p$-periodic solution (so that  the matrix $I-X(p)$ is singular). Let the vector $b$ be defined by (\ref{5.1}).

\noindent
Case 1. $b$ does not belong to the range of $I-X(p)$. Then all solutions of (\ref{3}) are unbounded as $t \ra \infty$. (The classical Massera's Theorem \ref{thm:0}.)

\noindent
Case 2. $b$  belongs to the range of $I-X(p)$. Then  (\ref{3}) has  infinitely many $p$-periodic solutions. Further sub-cases are as follows.

\medskip

\noindent
(i) If moreover $\rho (X(p))>1$, then (\ref{3})  has also unbounded solutions.
\medskip

\noindent
(ii) Assume that $\rho (X(p))=1$, and $\la =1$ is the only eigenvalue of $X(p)$ on the unit circle $|\la|=1$, and it has as many linearly independent eigenvectors as its multiplicity (i.e., the the Jordan block corresponding to $\la=1$ is diagonal). Then every solution of   (\ref{3}) approaches one of its $p$-periodic solutions, as $t \ra \infty$.
\medskip

\noindent
(iii) Suppose that $\rho (X(p))=1$, and there are other  eigenvalues of $X(p)$ on the unit circle $|\la|=1$, in addition to $\la =1$. Assume that all  eigenvalues of $X(p)$ on the unit circle $|\la|=1$ have diagonal Jordan blocks. Then all solutions of   (\ref{3}) are bounded, as $t \ra \infty$. 
\end{thm}

\pf
Let $x(t)$ be any solution of  (\ref{3}), represented by (\ref{4.1}). We shall consider the iterates $x(mp)$, where $m$ is a positive integer. With $b$ as defined by (\ref{5.1})
\[
x(p)=X(p)x_0+b \,.
\]
By periodicity, $x(t+p)$ is also a solution  of  (\ref{3}), which is equal to $x(p)$ at $t=0$.  Using  (\ref{4.1}) again
\[
x(t+p)=X(t)x(p)+X(t) \int_0^t X^{-1}(s) f(s) \, ds \,.
\]
Then
\[
x(2p)=X(p)x(p)+b=X(p) \left(  X(p)x_0+b \right)+b=X^2(p)x_0+X(p)b+b \,.
\]
By induction, for any integer $m>0$,
\beq
\lbl{7}
x(mp)=X^m(p)x_0+\sum _{k=0}^{m-1} X^k(p) b \,.
\eeq

{\em Case} 1. Assume that  $b$ does not belong to the range of $I-X(p)$. Then the linear system (\ref{6}) has no solutions. Since $\det \left(I-X(p) \right)^T=\det \left(I-X(p) \right)=0$, it follows that the system
\beq
\lbl{8}
\left(I-X(p) \right)^Tv=0
\eeq
has non-trivial solutions, and we claim that it is possible to find a  non-trivial solution $v_0$ of (\ref{8}) for which the scalar product with $b$ satisfies  
\beq
\lbl{9}
(b,v_0) \ne 0 \,.
\eeq
Indeed, assuming otherwise,  $b$ would be  orthogonal to the null-space of $\left(I-X(p) \right)^T$, and then the linear system (\ref{6}) would be solvable by the Fredholm alternative, a contradiction.
From (\ref{8}), $v_0=X(p)^Tv_0$, then $X(p)^Tv_0=X^2(p)^Tv_0$, which gives
$v_0=X^2(p)^Tv_0$, and inductively we get 
\beq
\lbl{10}
v_0=X^k(p)^Tv_0, \s \mbox{for all positive integers $k$}.
\eeq
Then by (\ref{7})
\[
\left( x(mp),v_0 \right)=\left( X^m(p)x_0,v_0 \right)+\sum _{k=0}^{m-1} ( X^k(p) b,v_0)
\]
\[
=(x_0,X^m(p)^T v_0)+\sum _{k=0}^{m-1} (  b,X^k(p)^T v_0)=(x_0, v_0)+m(b, v_0) \ra \infty  \,,
\]
as $m \ra \infty$, in view of (\ref{9}).
\medskip

{\em Case 2}. Assume now that  $b$  belongs to the range of $I-X(p)$. Then the linear system (\ref{6}) has a solution denoted by $\bar x_0$, and $x(t,\bar x_0)$ is a $p$-periodic solution of (\ref{3}). Adding to it non-trivial solutions of the corresponding homogeneous system (\ref{4}) produces infinitely many $p$-periodic solutions
of (\ref{3}). 
\medskip

Turning to the sub-cases, from (\ref{6})
\beq
\lbl{11}
\bar x_0=X(p) \bar x_0+b \,.
\eeq
Then
\[
\bar x_0=X(p) \left( X(p) \bar x_0+b  \right)+b=X^2(p) \bar x_0+X(p)b +b\,.
\]
Continuing  to use the latest expression for $\bar x_0$ in (\ref{11}), obtain inductively
\[
\bar x_0=X^m(p)\bar x_0+\sum _{k=0}^{m-1} X^k(p) b \,,
\]
so that $\sum _{k=0}^{m-1} X^k(p) b=\bar x_0-X^m(p)\bar x_0$. Using this in (\ref{7}), obtain
\beq
\lbl{11.1}
x(mp)=\bar x_0+X^m(p)\left(x_0-\bar x_0 \right) \,.
\eeq
In case $\rho (X(p))>1$ (the sub-case (i)), we can choose a vector $x_0$ to make $x(mp)$ unbounded, producing an unbounded solution of (\ref{3}) (choose  $x_0-\bar x_0$ to  be an eigenvector of $X(p)$ corresponding to an eigenvalue $\la$, with $|\la|>1$). 
\medskip

In the sub-case (ii),  assume for simplicity that $X(p)$ has a complete set of eigenvectors $z_1,z_2, \dots, z_k,\ldots, z_n$, with $z_1,z_2, \dots, z_k$ corresponding to the eigenvalue $\la =1$ of multiplicity $k<n$, and the other eigenvectors corresponding to the eigenvalues with $|\la|<1$. Decomposing $x_0-\bar x_0=\sum _{i=1}^n c_iz_i$, obtain (since $|\la|<1$ for all eigenvalues other than $\la =1$)
\[
X^m(p)\left(x_0-\bar x_0 \right) \ra \sum _{i=1}^k c_iz_i \equiv y \,,
\]
where $y$ is an eigenvector of $X(p)$ corresponding to the eigenvalue $\la =1$. 
It follows by (\ref{11.1}) that for any $x_0$, $x(mp,x_0) \ra \bar x_0+y$, and then $x(t, mp+x_0) \ra x(t,\bar x_0+y )$, or $x(t+ mp,x_0) \ra x(t,\bar x_0+y )$, where 
 $ x(t,\bar x_0+y )$ is one of the $p$-periodic solutions of (\ref{3}).
For the general case, one uses the Jordan normal form of $X(p)$, replacing the eigenvectors corresponding to $|\la|<1$ with the generalized eigenvectors.
\medskip

In the sub-case (iii),  similar arguments show that the sequence $\{ x(mp \}$ is bounded for any solution $x(t)$ of (\ref{3}). We claim that then $x(t)$ is bounded. Indeed, solutions of (\ref{3}) can have only a limited change over one period, by continuity, so that an unbounded solution cannot have  the sequence  $\{ x(mp \}$ bounded.
\epf

The assumption of Theorem \ref{thm:1} that the homogeneous system (\ref{4}) has a $p$-periodic solution can be seen as a {\em  case of resonance}. The complementary case when (\ref{4}) does not have a $p$-periodic solution is easy. Then the matrix $I-X(p)$ is non-singular, and hence the non-homogeneous system (\ref{3}) has a unique $p$-periodic solution for any $f(t)$. The difference of any two solutions of (\ref{4}) satisfies (\ref{3}), and therefore this  $p$-periodic solution is stable if $\rho \left(X(p) \right)<1$, and unstable if $\rho \left(X(p) \right)>1$.

\section{Instability for a class of first order equations}
\setcounter{equation}{0}
\setcounter{thm}{0}
\setcounter{lma}{0}

We now consider  nonlinear perturbations of first order equations
\beq
\lbl{15}
x'+a(t)x+g(x)=f(t) \,,
\eeq 
with $g(x) \in C^1(R)$, and  $a(t), f(t) \in C(R)$, satisfying $a(t+p)=a(t)$ and $f(t+p)=f(t)$ for all $t$, and some $p>0$. We assume that
\beq \lbl{15.1}
\int_0^p a(t) \, dt=0 \,,
\eeq
so that the linear part of this equation is at resonance.
Again, we denote $\mu (t)=e^{\int _0^t a(s) \, ds}$, which by (\ref{15.1}) is a $p$-periodic function. The nonlinear term $g(x)$ is assumed to satisfy a condition of E.M. Landesman and A.C. Lazer \cite{L}: the limits $g(\infty)$ and $g(-\infty)$ exist and  
\beq
\lbl{16}
g(-\infty) <g(x)<g(\infty),  \s \mbox{for all $x \in (-\infty,\infty)$} \,.
\eeq

\begin{thm}\lbl{thm:2}
Assume that (\ref{15.1}) and (\ref{16}) hold. The equation (\ref{15}) has a $p$-periodic solution if and only if $f(t)$ satisfies
\beq
\lbl{17}
g(-\infty) \int_0^p \mu (t) \, dt<\int_0^p \mu (t) f(t) \, dt<g(\infty) \int_0^p \mu (t) \, dt \,.
\eeq
If in addition to (\ref{15.1}), (\ref{16}) and (\ref{17})  
\beq
\lbl{18}
g'(x)>0 \,, \s \mbox{for all $x \in R$} \,,
\eeq
then the equation (\ref{15}) has a unique $p$-periodic solution that attracts all other solutions as $t \ra \infty$.
\medskip

If the conditions (\ref{15.1}) and (\ref{16}) hold, but  (\ref{17}) fails, then all of the solutions of (\ref{15}) are unbounded as $t \ra \infty$, and as $t \ra -\infty$.
\end{thm}

\pf
Let $x(t)$ be a  $p$-periodic solution of (\ref{15}). Multiply  (\ref{15}) by the $p$-periodic  $\mu (t)>0$, then integrate over $(0,p)$. Integrate by parts, using  that $\mu'=\mu a(t)$ and $\mu (0)=\mu (p)=1$, to obtain
\beq
\lbl{18.1}
\int_0^p g(x) \mu(t) \, dt=\int_0^p f(t)\mu(t) \, dt \,.
\eeq
It follows that (\ref{17}) holds, in view of (\ref{16}).
\medskip

Conversely, assume  that (\ref{17}) holds. The existence of $p$-periodic solution of (\ref{15}) will follow by a simple fixed point argument. Indeed, write  solutions    of  (\ref{15}) as
\[
x(t)=\frac{1}{\mu (t)} x(0)+\frac{1}{\mu (t)}\left[ \int_0^t \mu (s) f(s) \, ds -\int_0^t g(x(s)) \mu(s) \, ds \right] \,.
\]
Observe that $a(t)$ is bounded from above and from below by continuity, and the same is true for $\mu (t)$ by periodicity.  Hence if $A>0$ is large, then $x(t,A)$ is large for all $t \in (0,p)$. Then $g(x(t))$ is close to $g(\infty)$, and the term in the square bracket is negative, and hence $x(p,A)<A$. Similarly, $x(p,-A)>-A$, for $A>0$ large. It follows that the continuous  Poincar\'{e} map $x_0 \ra x(p,x_0)$ takes the interval $(-A,A)$ into itself. There exists a fixed point, leading to a  $p$-periodic solution.
\medskip

Assume now that  the condition (\ref{17}) fails. Assume for definiteness that 
\beq
\lbl{20}
\int_0^p \mu (t) f(t) \, dt \geq g(\infty) \int_0^p \mu (t) \, dt \,,
\eeq
and the case when $\int_0^p \mu (t) f(t) \, dt \leq g(-\infty) \int_0^p \mu (t) \, dt$ is similar.
Let $x(t)$ be any solution of (\ref{15}). Multiply  (\ref{15}) by the $p$-periodic  $\mu (t)>0$, then integrate over $(0,p)$. Since $x(t)$ is no longer assumed to be periodic, integration by parts produces two extra terms. Similarly to (\ref{18.1}) obtain
\beqa 
\lbl{20.1}
& x(p)-x(0)=\int_0^p f(t)\mu(t) \, dt-\int_0^p g(x(t)) \mu(t) \, dt  \\ \nonumber
&>\int_0^p f(t)\mu(t) \, dt-g(\infty) \int_0^p \mu (t) \, dt  \equiv \al \geq 0 \,. \nonumber
\eeqa
Assume first that $\al >0$, i.e., the inequality in (\ref{20}) is strict. Then
\[
x(p)-x(0)>\al>0 \,.
\]
Apply a similar argument on $[p,2p]$, and use the periodicity of $\mu (t)$ and $f(t)$ to get
\[
x(2p)-x(p) >\al >0\,.
\]
so that $x(2p)-x(0)>2\al$. Then $x(mp)-x(0)>m\al$ for any integer $m>0$, and hence $x(t)$ is unbounded. In case $\al =0$, we have $x(p)-x(0)>0$ from (\ref{20.1}), so  that the  Poincar\'{e} map $x(0) \ra x(p,x(0))$ satisfies $ x(p,x(0))>x(0)$ for all $x(0) \in R$. The increasing sequence $\{x(mp) \}$ has to go to infinity, since otherwise it would have to converge to a limit, which is a fixed point of the  Poincar\'{e} map. But fixed points are not possible for a  map that takes any number into a larger one.
\medskip

Assume finally  that (\ref{15.1}), (\ref{16}), (\ref{17}) and (\ref{18}) hold. By above, there is a $p$-periodic solution of (\ref{15}), call it $y(t)$.  Let $x(t)$ be any other solution of (\ref{15}), and set  $z(t)=x(t)-y(t)$. By the mean value theorem $z(t)$ satisfies a linear equation
\[
z'+b(t)z =0\,,
\]
with $p$-periodic $b(t)=a(t)+\int_0^1 g' \left( sx(t)+(1-s)y(t) \right) \, ds>a(t)$, so that 
$\int _0^p b(s) \, ds>0$. It follows that $z(t) \ra 0$, as $t \ra \infty$. In particular, this implies that  the periodic solution $y(t)$ is unique, and it attracts all other solutions as $t \ra \infty$.
\epf  

The argument we gave above for the  case $\al =0$ could be replaced by using the following result of J.M. Alonso and R. Ortega \cite{O} (the Corollary 2.3 in \cite{O}).

\begin{prop} (\cite{O}) \lbl{prop:o}
Consider a difference equation on a finite dimensional Banach space $X$:
\[
\xi _{n+1}=F \left(\xi _{n} \right), \s\s n \geq 0 \,,
\]
where $F: X \ra X$ is a continuous operator. If there exists a continuous functional $V$ satisfying
\[
V \left( F(\xi)\right)>V(\xi) \,, \s\s \forall  \xi \in X \,,
\]
then $\lim _{n \ra \infty} ||\xi _{n}||=\infty$.
\end{prop}

\noindent
{\bf \large Example } Consider an equation with the linear part at resonance
\beq
\lbl{21}
x'(t)+\sin t \, x(t)+ \frac{2}{\pi} \tan ^{-1} x(t)=\nu +\sin t \,,
\eeq
where $\nu$ is a parameter. Here $p=2 \pi$, $a(t)=\sin t$, $f(t)=\nu +\sin t$, $g(x)=\frac{2}{\pi} \tan ^{-1} x$, so that $g(-\infty)=-1$ and $g(\infty)=1$, with $g'(x)>0$. Calculate $\mu (t)=e^{1-\cos t}$, $\int _0^{2 \pi} \mu (t) f(t) \, dt=\nu \int _0^{2 \pi}  \mu (t) \, dt$. The condition (\ref{17}) becomes 
\[
-1<\nu<1 \,.
\]
Theorem \ref{thm:2} leads to the following conclusion: If $\nu \in (-1,1)$ the equation (\ref{21}) has a unique $2 \pi$-periodic solution that attracts all of its other solutions as $t \ra \infty$. If $\nu \geq 1$ or $\nu \leq -1$, then all solutions of  (\ref{21}) are unbounded, both  as $t \ra \infty$ and   as $t \ra -\infty$. 
\medskip

\begin{figure}
\begin{center}
\scalebox{0.9}{\includegraphics{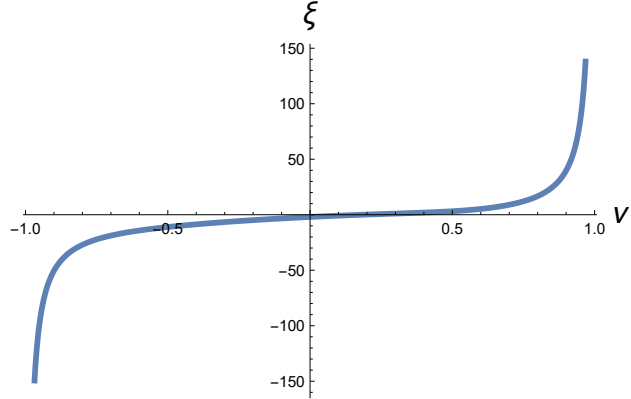}}
\caption{ The  curve of $2\pi$ periodic solutions of (\ref{21}), with their averages $\xi$ drawn versus $\nu$}
\label{fig:1}
\end{center}
\end{figure}

It turns out that  $2\pi$-periodic solutions of (\ref{21}) tend to infinity as $\nu \ra \pm 1$, see the Figure \ref{fig:1}. In that figure $\xi$ is the average of $2\pi$-periodic solutions $x(t)$, so that $x(t)=\xi+X(t)$, with $\int_0^{2\pi} X(t) \, dt=0$ (see the Theorem \ref{thm:last} below.) We used a modification of the {\em Mathematica} program presented and explained in \cite{KS1}.
\medskip

We thus obtained an exhaustive description of the dynamics of (\ref{21}), easily confirmed by numerical experiments.

\medskip

\noindent
{\bf \large Remark $\;$} Suppose that the condition (\ref{17}) holds, but (\ref{18}) does not. Then the equation (\ref{15}) has a $p$-periodic solution, but the asymptotic behavior of other solutions is an open problem.

\section{Unbounded solutions for a class of systems}
\setcounter{equation}{0}
\setcounter{thm}{0}
\setcounter{lma}{0}
\setcounter{prop}{0}

We recall some basic results on linear periodic systems. Consider the {\em adjoint system} for the  homogeneous $p$-periodic system (\ref{4}) 
\beq
\lbl{4ad}
z'=-A^T(t)z \,,
\eeq
where $A^T$ denotes the transpose. The following two lemmas can be found in B.P. Demidovi\v{c} \cite{d}. We include slightly simpler proofs for completeness.

\begin{lma}\lbl{lma:ad}
If the  system (\ref{4}) has a non-trivial $p$-periodic solution, then so does (\ref{4ad}).
\end{lma}

\pf
Let $X(t)$ be again the fundamental solution matrix of (\ref{4}). We are given that $X(p)$ has an eigenvalue $\la =1$. Recall that
\beq
\lbl{90}
X'=A(t)X \,.
\eeq 
Let $Z(t)$ be the fundamental solution matrix of (\ref{4ad}), so that  
\beq
\lbl{91}
Z'=-A^T(t)Z \,.
\eeq
We claim that $Z=Y^{-1}(t)$, where $Y(t)=X^T(t)$, i.e., $Z=\left(  X^T \right)^{-1}$. Indeed, using that $Z'=-Y^{-1}Y'Y^{-1}$ (differentiate $YY^{-1}=I$, or see e.g., p. 5 in R. Bellman \cite{be}), in order to justify that $Z$ satisfies (\ref{91}) the following equivalent statements must hold:
\[
-Y^{-1}Y'Y^{-1}=-A^T Y^{-1} \,,
\]
\[
-Y^{-1}Y'=-A^T  \,,
\]
\[
Y'=YA^T \,,
\]
\[
\left(  X^T \right)'=X^T A^T  \,,
\]
\[
X'=AX \,,
\]
which is (\ref{90}), proving the claim. The eigenvalues of $Y(p)$ are the same as  those of $X(p)$, so that one of them is $\la =1$. The eigenvalues of $Z(p)$ are the reciprocals of those of $Y(p)$, so that one of them is $\la =1$, and the  system (\ref{4ad}) has a $p$-periodic solution.
\epf

\begin{lma}\lbl{lma:ad11}
Assume that the homogeneous system (\ref{4})  has a $p$-periodic solution. Then the non-homogeneous system (\ref{3})  has a $p$-periodic solution if and only if the integral of scalar product
\beq
\lbl{30}
\int _0^p f(t) \cdot z(t) \, dt=0 \,,
\eeq
where $z(t)$ is any  $p$-periodic solution of (\ref{4ad}).
\end{lma}

\pf
Let $x(t)$ and $z(t)$  be $p$-periodic solutions  of (\ref{3}) and (\ref{4ad}) respectively. 
To prove the necessity part, multiply the equation $i$ of (\ref{3}) by $z_i$,  the equation $i$ of (\ref{4ad}) by $x_i$, add the equations, and sum in $i$. Then integrate, and use the periodicity to obtain:
\[
\int _0^p f(t) \cdot z(t) \, dt=\int _0^p \left[Ax \cdot z-  x \cdot A^T z \right]  \, dt =0 \,.
\]

Turning to the sufficiency part,  any  $p$-periodic solution of (\ref{4ad}) can be written as $z(t)=Z(t)z_0$, where $z_0$ satisfies
\[
\left[ I-Z(p) \right]z_0=\left[ I-\left( {X^T} \right)^{-1}(p)  \right]z_0=0 \,,
\]
which can be written as
\beq
\lbl{31}
z_0^TX(p)=z_0^T \,,
\eeq
or as 
\beq
\lbl{32}
\left( I-X^T(p)   \right)z_0=0 \,.
\eeq
We are given that (\ref{30}) holds, which can be written as
\beq
\lbl{33}
0=\int _0^p   \left(   Z(t)z_0   \right)^T f(t) \, dt=z_0^T \int _0^p   X^{-1}(t) f(t) \, dt \,.
\eeq
Since the system (\ref{4})  has a $p$-periodic solution, both (\ref{5}) and (\ref{32}) have non-trivial solutions. In order for  (\ref{3})  to have a $p$-periodic solution, the system of equations (\ref{6}) has to be solvable, which requires that the vector $b$ defined in (\ref{5.1}) must be orthogonal to any solution $z_0$ of (\ref{32}). Using (\ref{31}) and (\ref{33}), obtain
\[
b \cdot z_0=z_0^Tb=z_0^T X(p) \int _0^p   X^{-1}(t) f(t) \, dt=z_0^T  \int _0^p   X^{-1}(t) f(t) \, dt=0 \,.
\]
completing the proof.
\epf

We wish to extend the Theorem \ref{thm:2} to systems. This can be done in a number of ways. For example, a recent paper of A.  Boscaggin et al \cite{bo} considered coupled harmonic oscillators, each one at resonance. They show that solutions are unbounded if a condition of   A.C. Lazer  and  D.E. Leach \cite{leach} type is violated. We shall obtain a straightforward extension of the Theorem \ref{thm:2} (which used a condition of  E.M. Landesman and A.C. Lazer \cite{L}) provided that the adjoint system (\ref{4ad}) has {\em a positive} $p$-periodic solution, and give a condition for that to happen.
\medskip

We consider bounded nonlinear perturbations of linear systems
\beq
\lbl{40}
x'+A(t)x+f(x)=g(t) \,.
\eeq
Here   an $n \times n$ matrix $A(t)$, and $g(t) \in R^n$  have  continuous $p$-periodic entries, the unknown vector $x=x(t)  \in R^n$,  
\[
f(x)=
\left[
\begin{array}{c}
 f_1(x) \\
 \vdots \\
 f_n(x)
\end{array}
\right]
\]
is a continuous vector function.
We assume that the linear part is at resonance, so that both 
\beq
\lbl{41}
x'+A(t)x=0 \,,
\eeq
and 
\beq
\lbl{42}
z'-A^T(t)z=0
\eeq
have non-trivial $p$-periodic solutions, and moreover that $z(t)>0$ componentwise ($z_i(t)>0$ for all $i$). The components of the vector $f(x)$ are assumed to satisfy
\beq
\lbl{43}
\al _i <f_i(x)<\beta _i \,, \s\s \mbox{for all $x \in R^n$, and all $i$} \,,
\eeq
with $2n$ given constants $\al _i, \beta _i$.

\begin{thm}\lbl{thm:3}
Assume that the adjoint system (\ref{42}) has a positive $p$-periodic solution $z(t)$, and (\ref{43}) holds. Then the system (\ref{40}) may have a $p$-periodic solution only if 
\beq
\lbl{44}
\sum _{i=1}^n \al _i \int _0^p z_i(t) \, dt<\int _0^p g(t) \cdot z(t) \, dt< \sum _{i=1}^n \beta _i \int _0^p z_i(t) \, dt \,.
\eeq
In case this condition fails, then all solutions of  (\ref{40}) are unbounded as $t \ra \pm \infty$.
\end{thm}

\pf
If $x(t)$ is a  $p$-periodic solution of  (\ref{40}), then $f(x(t))$ is a  $p$-periodic function. Applying Lemma \ref{lma:ad11} obtain
\[
\int _0^p f(x(t)) \cdot z(t) \, dt=\int _0^p g(t) \cdot z(t) \, dt \,,
\]
from which (\ref{44}) follows, since $z(t)>0$.
\medskip

Assume now that the condition (\ref{44}) fails. Suppose for definiteness that
\[
\sum _{i=1}^n \al _i \int _0^p z_i(t) \, dt \geq \int _0^p g(t) \cdot z(t) \, dt \,.
\]
Multiply the equation $i$ in (\ref{40}) by $z_i(t)$, integrate over $(0,p)$ then add up in $i$. Integrating by parts, and using $p$-periodicity of $z(t)$ and (\ref{43}), obtain
\[
\sum _{i=1}^n z_i(0) \left[ x_i(p)-x_i(0) \right]=\int _0^p g(t) \cdot z(t) \, dt-\int _0^p f(x) \cdot z(t) \, dt 
\]
\[
  <\int _0^p g(t) \cdot z(t) \, dt-\sum _{i=1}^n \al _i \int _0^p z_i(t) \, dt \leq 0 \,.
\]
So that
\[
-\sum _{i=1}^n z_i(0)  x_i(p) >-\sum _{i=1}^n z_i(0) x_i(0)   \,.
\]
We now apply the Proposition \ref{prop:o}, with $X=R^n$, the Poincar\'{e} map $F: x(0) \ra x(p)$, and the functional $V\left( x(t) \right) =-\sum _{i=1}^n z_i(0) x_i(t)$ to conclude the unboundness of the sequence $\{x(mp) \}$.
\epf

\noindent
{\bf Remark } We do not know if the conditions (\ref{43}) and (\ref{44}) are sufficient for the existence of $p$-periodic solutions of (\ref{40}). As mentioned in A.  Boscaggin  et al \cite{bo}, few existence results are known for semilinear periodic systems.
\medskip

To give a condition for (\ref{42}) to have a positive $p$-periodic solution we need the following lemma.

\begin{lma}\lbl{lma:pos}
Assume that  $B(t)$ is a continuous $n \times n$ matrix with positive off diagonal entries for all $t>0$, and $e_i \in R^n$ has the entry $i$ equal to one, and the other entries are zero. Then solution of
\beq
\lbl{45}
y'=B(t)y \,, \s\s y(0)=e_i
\eeq
satisfies $y(t)>0$ for all $t>0$, and all $i=1,2,\ldots,n$.
\end{lma}

\pf
We can find a constant matrix $B_0$ with  positive  off diagonal entries, such that $B(t)>B_0$ for small $t$. It is well known that solutions of
\[
x'=B_0x \,, \s\s x(0)=e_i
\]
are positive for all $t>0$, see e.g., p. 176 in R. Bellman \cite{be}. Since $y(t)>x(t)$ componentwise, it follows that
  $y(t)>0$ for small $t$. At the first  $t_0$ where $y_k(t_0)=0$ for some $k$,  there is a contradiction in the $k$-th equation  of (\ref{45}) at $t=t_0$, since $y_k'(t_0)>0$ from (\ref{45}).
Hence, $y(t)>0$ for all $t>0$.
\epf

\begin{prop}
Assume that  off diagonal entries of a $p$-periodic matrix $A(t)$ are positive, and the spectral radius $\rho (Z(p))=1$. Then (\ref{42})  has a positive $p$-periodic solution.
\end{prop}

\pf
By Lemma \ref{lma:pos} the fundamental matrix $Z(p)$ of  (\ref{42})  has  positive entries. By the Perron-Frobenius theorem the largest in absolute value eigenvalue of $Z(p)$ is positive, and since $\rho (Z(p))=1$, it is $\la =1$, and the corresponding eigenvector $\xi$ is also positive. Then $Z(t) \xi$ gives positive $p$-periodic solution of  (\ref{42}).
\epf

For the  $2 \times 2$ case we give  conditions that appear easier to check.

\begin{prop}
Assume that  $A(t)$ is  a $p$-periodic $2 \times 2$ matrix with $a_{12}(t)>0$ and  $a_{21}(t)>0$ for all $t$, and $\int _0^p \left[a_{11}(t)+a_{22}(t)   \right] \, dt \leq 0$. Finally, assume that $Z(p)$ has an eigenvalue $\la =1$.  Then (\ref{42})  has a positive $p$-periodic solution.
\end{prop}

\pf
By Lemma \ref{lma:pos} the fundamental matrix $Z(p)$ of  (\ref{42})  has  positive entries. If $\la _1$ and $\la _2$ are the eigenvalues of $Z(p)$, then $\la _1=1$, and by Liouville's formula, see e.g., p. 212 in \cite{K},
\[
0<{\rm Det } \, X(p) =\la _1 \la _2=e^{\int _0^p \left[a_{11}(t)+a_{22}(t)   \right] \, dt} \leq 1,
\] 
so that $0 \leq \la _2 \leq 1$. By the Perron-Frobenius theorem $0 \leq \la _2 < 1$, and the  eigenvector $\xi$ corresponding to $\la _1 =1$ is  positive. Then $Z(t) \xi$ gives a positive $p$-periodic solution of  (\ref{42}).
\epf

\section{Solution curves and unboundness of solutions}
\setcounter{equation}{0}
\setcounter{thm}{0}
\setcounter{lma}{0}

We now consider  nonlinear perturbations of a second order  periodic problem at resonance
\beq
\lbl{60}
x''(t)+\la x'(t)+g(x)=f(t) \,, 
\eeq 
with $g(x) \in C^1(R)$, and  $a(t), f(t) \in C(R)$, satisfying $a(t+p)=a(t)$ and $f(t+p)=f(t)$ for all $t$ and some $p>0$, and a constant $\la >0$. This pendulum-like equation was studied previously in a number of papers, including  J. \v{C}epi\v{c}ka et al \cite{ce}, G. Tarantello \cite{ta}, A. Castro \cite{ca}. As mentioned above, the linear part of this equation (when $g(x) \equiv 0$) is at resonance. Decompose $f(t)=\mu +e(t)$, with $\mu \in R$ and $\int_0^p e(t) \, dt=0$. Similarly, decompose the solution  $x(t)=\xi +X(t)$, with $\xi \in R$ and $\int_0^p X(t) \, dt=0$.  In view of the above decomposition, we may write (\ref{60}) as 
\beq
\lbl{61}
x''+\la x'+g(x)=\mu +e(t) \,.
\eeq

The following result we proved in \cite{K2}.
\begin{thm}\lbl{thm:5}
Assume that $g(x) \in C^1(R)$ is a bounded function ($|g(x)| \leq M$ for all $x \in R$ and some $M>0$), and 
\beq
\lbl{62}
|g'(x)|<\frac{\la ^2}{4} +\omega ^2 \,, \s \mbox{for all $x \in R$} \,, \s\s \mbox{where $\omega =\frac{2 \pi}{p}$} \,.
\eeq
Then for any $\xi \in R$ one can find a unique $\mu \in R$ for which the problem (\ref{61}) has a unique $p$-periodic solution. Moreover, all  $p$-periodic solutions of  (\ref{61}) lie on a unique continuous solution curve $(\mu ,x(t))(\xi)$.
\end{thm}

We now give an instability result based on the Landesman-Lazer \cite{L} condition.
\begin{thm}\lbl{thm:curva}
In addition to the conditions of the Theorem \ref{thm:5} assume that the limits at infinity $g(\pm \infty)$ exist, and
\beq
\lbl{63}
g(-\infty)<g(x)<g(\infty) \,, \s \mbox{for all $x \in R$}.
\eeq 
Then the equation (\ref{61}) has a $p$-periodic solution if and only if 
\beq
\lbl{64}
g(-\infty)< \mu <g(\infty)  \,.
\eeq
If the condition (\ref{64}) fails, then all of the solutions of (\ref{61}) are unbounded as $t \ra \infty$ and as $t \ra -\infty$.
\end{thm}

\pf
Let $x(t)$ be  a $p$-periodic solution of (\ref{63}). 
Integrate the equation (\ref{61}) over $(0,p)$:
\beq
\lbl{65}
\mu p=\int_0^p g(x(t)) \, dt \,.
\eeq
Then the necessity of the condition (\ref{64}) follows by (\ref{63}).
\medskip

By the Theorem \ref{thm:5} there is a continuous solution curve $(\mu ,x(t))(\xi)$ for $\xi \in R$. Moreover, we showed in \cite{K2} that with $x(t)=\xi+X(t)$, there is a uniform in $\xi$ and $t$ bound on $|X(t)|$. It follows from (\ref{65}) that $\mu \ra g(\infty)$ ($\mu \ra g(-\infty)$) as $\xi \ra \infty$ ($\xi \ra -\infty$). By the continuity of the solution curve, it follows that  the condition (\ref{64})  is sufficient for the existence of $p$-periodic solution.
\medskip

If the condition (\ref{64}) fails, assume for definiteness that
\beq
\lbl{66}
\mu \geq g(\infty)  \,.
\eeq
If $x(t)$ is any solution of (\ref{61}), integration of this equation gives
\beqa \nonumber
& x'(p)-x'(0)+\la \left(  x(p)-x(0)  \right)=\mu p-\int_0^p g(x(t)) \, dt  \\ \nonumber
& > \mu p-g(\infty) p  \geq 0 \,, \nonumber
\eeqa
so that $x'(p)+\la x(p)>x'(0)+\la x(0)$.
We now apply the Proposition \ref{prop:o}, with $X=R^2$, the Poincar\'{e} map $F: \left( x(0),x'(0) \right) \ra \left( x(p),x'(p) \right)$, and the functional $V\left( x(t),x'(t) \right) =x'(t)+\la x(t)$ to conclude the unboundness of the sequence $\{x(mp),x'(mp) \}$, proving  the unboundness of solutions of (\ref{61}).
\epf

Similar results hold for  first order periodic equations of the type
\beq
\lbl{la}
x'(t)+g(x)=\mu +e(t) \,, 
\eeq
with $\mu \in R$, $e(t) \in C(R)$, satisfying $e(t+p)=e(t)$ for all $t$ and some $p>0$, and $\int _0^p e(t) \, dt=0$. As above, decompose the $p$-periodic solutions of (\ref{la}) as  $x(t)=\xi +X(t)$, with $\xi \in R$ and $\int_0^p X(t) \, dt=0$. We shall sharpen our result in \cite{K2} on the solution curve $(\mu,x(t))(\xi )$, after establishing the following lemma.
\begin{lma} \lbl{lma:1a}
Consider a linear periodic problem in the class of functions of zero average
\beq
\lbl{lb}
w'(t)+h(t)w(t)=\mu \,, \s w(t+p)=w(t)\,, \s   \s \int _0^p w(t) \, dt=0 \,,
\eeq
where $h(t) \in C(R)$ is a given function of period $p$, and $\mu$ is a parameter. The only solution of (\ref{lb}) is $\mu=0$ and $w(t) \equiv 0$.
\end{lma}

\pf
We claim that $w(t)$ is of one sign. If $\mu=0$ this follows by the explicit solution. If, say, $\mu >0$ and $w(t)$ is a sign changing solution, then by the periodicity of $w(t)$ one can find a point $t_0$ such that $w(t_0)=0$ and $w'(t_0) \leq 0$, which contradicts the equation (\ref{lb}). Since  $w(t)$ is of one sign and  of zero average, $w(t) \equiv 0$, and then  $\mu=0$ from  the equation (\ref{lb}).
\epf

\begin{thm}\label{thm:last}
Assume that $g(x) \in C^1(R)$, $e(t) \in C(R)$ is $p$-periodic of  zero average. Then for any $\xi \in R$ one can find a unique $\mu \in R$ for which the problem (\ref{la}) has a unique $p$-periodic solution. Moreover, all  $p$-periodic solutions of  (\ref{la}) lie on a unique continuous solution curve $(\mu ,x(t))(\xi)$.
\end{thm}

\pf
Local properties of the solution curve, and the fact that $\xi$ is a global parameter, were proved in \cite{K2}. We show next that $\mu$ and $x(t)$ are bounded, when $\xi$ belongs to a bounded set, so that the solution curve can be continued globally, for $-\infty<\xi<\infty$. With  $x(t)=\xi +X(t)$, obtain
\beq
\lbl{lc}
\s\s \s\s X'(t)+g(\xi+X(t))=\mu +e(t) \,, \s X(t+p)=X(t) \,,  \s \int _0^p X(t) \, dt=0 \,.
\eeq
Multiply the equation in (\ref{lc}) by $X'$ and integrate over $(0,p)$. By periodicity of $X(t)$
\[
\int _0^p {X'}^2(t) \, dt=\int _0^p X'(t) e(t) \, dt \,,
\]
which gives a bound on $\int _0^p {X'}^2(t) \, dt$. By Wirtinger's inequality obtain a bound on $\int _0^p X^2(t) \, dt$. With $X(t)$ bounded in $H^1$ norm, conclude a uniform bound on $X(t)$, and hence on $x(t)$, by Sobolev embedding. From (\ref{la})
\[
\mu ^2 \leq c_0 \left( {X'}^2(t)+g^2(x(t))+e^2(t) \right) \,,
\]
with some $c_0>0$. Integrating over $(0,p)$, obtain a bound on $\mu$.
\epf

For the equation (\ref{21}), considered above, we computed the $\mu =\mu (\xi)$ section of the solution curve described in Theorem \ref{thm:last}, and then plotted the inverse function $\xi=\xi (\mu)$ to produce the Figure \ref{fig:1}. The Theorem \ref{thm:last} can also be used to provide an alternative proof of the Theorem \ref{thm:2} (with extra information on the solution curve), similarly to the Theorem \ref{thm:curva}.

\end{document}